\documentclass[a4paper,11pt]{article}
\usepackage[latin2]{inputenc}
\usepackage[mathscr]{eucal}
\usepackage{t1enc}
\usepackage{amssymb}
\usepackage{amsfonts}
\usepackage{amsmath}
\usepackage{amsthm}
\usepackage{latexsym}
\normalfont
\usepackage[dvips]{color}

\usepackage[dvips]{epsfig}

\numberwithin{equation}{section}

\newtheorem{theorem}{Theorem}
\newtheorem*{theoremquote}{Theorem}

\newtheorem{lemma}[theorem]{Lemma}

\newtheorem{corollary}[theorem]{Corollary}

\theoremstyle{definition}

\theoremstyle{definition}

\newcommand{\eqn}[2]{\begin{equation}\label{#1}#2\end{equation}}
\newcommand{\eqnst}[1]{\begin{equation*}#1\end{equation*}}
\newcommand{\eqnspl}[2]{\begin{equation}\begin{split}\label{#1}%
	#2\end{split}\end{equation}}
\newcommand{\eqnsplst}[1]{\begin{equation*}\begin{split}%
	#1\end{split}\end{equation*}}

\def\cC{{\cal C}}
\def\hR{{\hat R}}
\def\hV{{\hat V}}

\def\Z{{\mathbb Z}}
\def\Ztwo{{{\mathbb Z}^2}}
\def\cnctd{\longleftrightarrow}
\newcommand\pcnctd[1]{\stackrel{#1}{\cnctd}}
\def\shcnctd{\leftrightarrow}

\title{The size of a pond in 2D invasion percolation}
\author{
Jacob van den Berg\footnote{Research funded in part by the Dutch BSIK/BRICKS project.}, \, 
Antal A. J\'{a}rai\footnote{Reserach supported by NSERC of Canada.} \,
 and Bálint Vágvölgyi
\\[0.1in]
{\small CWI and VUA, Carleton University and VUA} \\
{\footnotesize email: J.van.den.Berg@cwi.nl; jarai@math.carleton.ca; bvagvol@few.vu.nl}
}
\date{August 29, 2007}

\begin{document}

\maketitle

\begin{abstract}
We consider invasion percolation on the square lattice.
In \cite{BPSV} it has been proved that the probability that the radius of a so-called
pond  is larger than $n$, differs at most a factor of order $\log n$ 
from the probability that in critical Bernoulli percolation the radius of
an open cluster is larger than $n$.
We show that these two probabilities are, in fact, of the 
same order. Moreover, we prove an analogous result for the volume of a 
pond. 
\end{abstract}

\section{\textbf{Introduction}}
\label{main}
\parindent=5mm
Invasion percolation is a stochastic growth model of an interesting self-organised critical nature:
it has characteristics that resemble critical Bernoulli percolation,
even though the definition of the invasion process does not
involve any parameter (see \cite{newman,wilkinson}).
Comparison of the 
two processes helps to gain new insights into both of them (see e.g.
\cite{newman2,jarai,LPS06,AGdHS06}). 

Recently a new comparison result,
relating a so-called `pond' in invasion percolation to a critical Bernoulli
percolation cluster, has been proved in \cite{BPSV}. This result is
sharpened and extended in the present paper.

In the remainder of this Section we define 
the invasion percolation model and state
our main results. The proofs, and important prerequisites, are given
in Section \ref{proof}. 

For general background on percolation, see \cite{grimmett}.

\smallskip\noindent
Consider the hypercubic lattice
$\mathbb{Z}^d$ with its set of nearest neighbour bonds
$\mathbb{E}^d$. If an edge $e$ has endpoints $v$ and $w$, we write 
$e=\langle v,w \rangle$.
For an arbitrary subgraph $G = (V,E)$ of $(\mathbb{Z}^d,\mathbb{E}^d)$,
we define the outer boundary $\Delta G$ as
\begin{eqnarray*}
  \Delta G
  = \{ e = \langle v,w\rangle \in \mathbb{E}^d\ :\ \text{$e \notin E$, but 
      $v \in V$ or $w \in V$} \}.
\end{eqnarray*}
Invasion percolation is defined as follows. Let $\tau(e), e \in \mathbb{E}^d$, be independent
random variables, uniformly distributed on the interval $[0,1]$.
Using these variables,
we construct inductively an increasing sequence 
$G_0, G_1, G_2,\ldots$ of connected subgraphs of
the lattice. $G_0$ only contains the origin.
If $G_i = (V_i,E_i)$ has already been defined, we select the bond 
$e_{i+1}$ which minimizes $\tau$ on $\Delta G_i$, take
$E_{i+1} = E_i \cup \{e_{i+1}\}$, and let $G_{i+1}$ be the graph induced by 
the edge set $E_{i+1}$. The graph $G_i$ is called the invaded cluster at time $i$, 
and $G_\infty = (V_\infty,E_\infty) = \cup_{i=0}^{\infty}G_i$ is the
invaded region at time infinity.

Invasion percolation can be coupled to
Bernoulli bond percolation in the following standard way. Let $0\leq
p\leq 1$. For each bond $e$ we say that $e$ is $p$-open, if
$\tau(e)<p$. One can then define, in an obvious way, $p$-open paths and $p$-open clusters, and the
study of these objects corresponds with Bernoulli bond
percolation with parameter $p$.

From now on we will only consider the case when $d=2$. It is well-known
and easy to see that for every $p\in[0,1]$ the following holds:
once the invasion reaches an infinite $p$-open cluster,
it never leaves it again.
Further, it is a classical result
for $2D$ Bernoulli percolation that
for every $p > p_c$ (which equals $1/2$ on the square lattice)  there is (a.s.) a $p$-open circuit that surrounds
$O$ and belongs to the infinite $p$-open cluster, and that
(a.s.) there is no infinite
$p_c$-open cluster.
These properties
easily imply that (a.s.) $\hat\tau := \max_{e\in E_\infty} \tau(e)$ exists 
and is larger than $p_c$. Let $\hat{e}$ denote the edge where the 
maximum is taken and suppose that it is added to the invasion cluster 
at step $\hat{i} + 1$. Following the terminology in \cite{St-New},
the graph $G_{\hat{i}} = (\hat{V}, \hat{E})$ 
is called a `pond', or, more precisely, the first pond of $O$.
Since the invasion can be started at any vertex
$v$, not necessarily $O$, we have the more general notion 
`first pond of $v$'.

The above defined `pond' is a very natural object 
(see \cite{St-New} and \cite{BPSV}), and has several 
interpretations, for instance the following.
In this (somewhat informal) interpretation each vertex $(x,y) \in \Z^2$ represents a `polder': the square piece of (flat) land
$(x - 1/2, x + 1/2) \times (y-1/2,y+1/2)$, surrounded by four dikes, corresponding with
(the dual edges of) the four
edges of $(x,y)$. The heights of the dikes are the $\tau$ values of the corresponding edges. Now suppose that
water is supplied from some external source to the polder represented by the vertex $O$.
The water in this polder
will rise until its level reaches the height of the lowest of its four dikes, say $a$.
Then the water starts spilling
over that dike, so that the level in the neighbouring polder (on the other side of the dike) starts to rise.
If each of the other three
dikes of that neighbouring polder is higher than $a$, the water in the polder
of $O$ will remain at level $a$ until the above mentioned neighboring polder has reached this same
water level, 
after which the level in both polders rises (`simultaneously') until it reaches the height of the
lowest of the six dikes
bounding the union of these two polders, etc. On the other hand, if the neighboring polder of $O$
has a dike with height $b < a$, the water level in this polder will rise up to level $b$ and then
starts spilling over that dike to a third polder (while the level in the polder of $O$ is still $a$) etc.
In any case, (a.s.) eventually the water level at $O$ will remain constant forever, namely at the
level $\hat\tau$ defined above,
and the `connected' set of polders with the same final level as $O$ is the above
defined `first pond'.
Since water keeps being supplied, the `surplus' water will spill over the lowest dike (corresponding
with the earlier defined $\hat e$), on the boundary of this pond: the outlet from this first pond to a second
(lower level) pond. For clarity we note that for each vertex in the latter pond,
this pond plays the role of `first pond'. 

For further clarity we also note that the above `hydrologic' interpretation has
a more `symmetric' version as follows: Now at {\it each} vertex there is an external water source (rain, e.g.).
Again each polder has a final water level, and the maximal connected set of polders with the same final
water level, containing a given vertex $v$, is the earlier defined (first) pond of $v$.
Then, if $\hat V(v)$ denotes the vertices of the first pond of $v$, the collection
$\{ \hat V(v) \}_{v \in \Ztwo}$ is a random partition of $\Ztwo$ which is stationary under translations.

Before stating the results, we first fix some notation. 
Let $\hat{R} := \max\{|x| + |y| \, \colon \, (x,y) \in \hat V\}$ be the
radius of the first pond. Let $P_{cr}$ denote the product measure corresponding to
critical Bernoulli bond percolation. 
Let $B(n)$ denote the box $[-n,n]^2$ and $\partial B(n) := B(n) \setminus B(n-1)$. \\
Let $A$ and $B$ be sets of vertices.
In the context of Bernoulli percolation, we denote the event that there is an open path from $A$ to $B$
by $\{A \leftrightarrow B\}$. In the context of invasion percolation we denote the event that there is
a $p$-open path from $A$ to $B$ by $\{A \stackrel p \leftrightarrow B \}$. To indicate that there is an
infinite open (or $p$-open) path from $A$, we use the same notation with $B$ replaced by $\infty$. \\
We use the notation 
$g(n)\approx h(n)$, $n\to\infty$ to indicate that
\begin{eqnarray*}
\frac{\log g(n)}{\log h(n)}\to 1,\quad\text{as $n \to\infty$},
\end{eqnarray*}
and $g(n) \asymp h(n)$ to indicate that $g(n)/h(n)$ is bounded away from $0$ an $\infty$.

\smallskip\noindent
Van den Berg, Peres, Sidoravicius and Vares have proved the following 
theorem:

\begin{theoremquote}
\textbf{\cite[Proposition 1.3]{BPSV}} 
\begin{eqnarray}
  P(\hat{R}\geq n)
  \approx P_{cr}(0\leftrightarrow\partial B(n)),\ n\to\infty.
\end{eqnarray}
\end{theoremquote}

Using ideas and techniques from \cite{jarai}, we obtain the 
following improvement of the theorem above.

\begin{theorem}
\label{thm:radius}
\begin{eqnarray}
P(\hat{R}\geq n) \asymp  P_{cr}(0\leftrightarrow\partial B(n)).
\end{eqnarray}
\end{theorem}

Moreover, we show that not only the radius but also the volume of the pond behaves like that of
a critical percolation cluster: Let
\eqnst
{ s(n) 
  = n^2 P_{cr} ( 0 \shcnctd \partial B(n) ) \qquad\qquad
  \cC(0) 
  = \{ v \in \Ztwo : 0 \shcnctd v \}. }

\begin{theorem}
\label{thm:volume}
There exist constants $0 < c, c' < \infty$, such that
\eqnspl{e:volume}
{ c P_{cr} ( 0 \shcnctd \partial B(n) )
  &\leq P_{cr} ( |\cC(0)| > s(n) ) 
  \leq P( |\hat{V}| > s(n)) \\
  &\leq c'  P_{cr} ( 0 \shcnctd \partial B(n) ). }
\end{theorem}

\begin{corollary}
\label{cor:volume}
\begin{equation} \label{corv}
P( |\hat{V}| \geq n) \asymp P_{cr} ( |\cC(0)| \geq n).
\end{equation}
\end{corollary}

{\bf Remark:} These results, and the proofs in Section \ref{proof} also hold
(with some obvious adaptations) for the triangular and the hexagonal lattice.

\section{Proofs of the main results}\label{proof}
In the following all the constants are strictly positive and finite
without further mentioning.
\subsection{Preliminaries}
Let
\begin{eqnarray*}
\sigma(n,m,p)=P(\text{there is a $p$-open horizontal crossing of $[0,n]\times[0,m]$}),
\end{eqnarray*}
where it is assumed that the crossing does not use bonds lying on
the top or the bottom sides of the rectangle. Given
$\varepsilon>0$, we define, for $p > p_c$,
\begin{eqnarray*}
L(p,\varepsilon)=\min\{n:\sigma(n,n,p)\geq 1-\varepsilon\}.
\end{eqnarray*}
It is shown in \cite[(1.24)]{kesten}, that there exists an
$\varepsilon_0>0$ such that for all $\varepsilon\leq\varepsilon_0$, 
the scaling of $L(p,\varepsilon)$ is independent of
$\varepsilon$ in the sense that for all fixed
$0<\varepsilon_1,\varepsilon_2\leq\varepsilon_0$ the ratio
$L(p,\varepsilon_1)/L(p,\varepsilon_2)$ is bounded away from both 0
and $\infty$ as $p\downarrow p_c$. We let
$L(p)=L(p,\varepsilon_0)$ for the entire proof. Below we list some
properties of $L(p)$ that will play a crucial role in the proof of our results.
The first two follow fairly easily from the definitions and standard arguments (see Section
2.2 in \cite{jarai} for further explanation and references). 
The third is (a consequence of) a deep result in \cite{kesten}
\begin{enumerate}
\item $L(p)$ is decreasing,
right continuous and $L(p)\to\infty$ as $p\downarrow p_c$.
\item 
There is a constant
$D$ such that
\begin{eqnarray}\label{smalljump}
\lim_{\delta\downarrow 0}\frac{L(p-\delta)}{L(p)}\leq D \quad
\forall\ p>p_c.
\end{eqnarray}
\item \textbf{Theorem \cite[Theorem 2]{kesten}} There are
constants $C_0 > 0$ and $C_1$ such that for all $p>p_c$ 
\begin{eqnarray}\label{thetacorrineq}
C_0 P_{cr}\big[0\leftrightarrow\partial B(L(p))\big] \leq \theta(p)
\leq C_1 P_{cr}\big[0\leftrightarrow\partial B(L(p))\big],
\end{eqnarray}
where $\theta(p)=P_p(0 \leftrightarrow \infty)$ is the
percolation function for Bernoulli percolation.
\end{enumerate}

Finally we mention the following result on the behavior of 
$P_{cr}(0\leftrightarrow \partial B(n))$.
It is believed (see Chapters 9 and 10 in \cite{grimmett} for background) that for $2D$ percolation
on sufficiently `nice' 2D lattices this has a power law
(with critical exponent 5/48) but so far this has only been proved for site percolation
on the triangular lattice (see \cite{LSW02}). The following is sufficient for our
purpose. \\
There exists a constant $D_1$ such that
\begin{eqnarray}
\label{csakegyszerkell}
  \frac{P_{cr}(0\leftrightarrow \partial B(n))}{P_{cr}(0\leftrightarrow \partial B(m))}
  \geq D_1\sqrt{\frac{m}{n}},\quad 1\leq m\leq n.
\end{eqnarray}
For $m=1$ this was proved in \cite[Corollary (3.15)]{BK}. For general $m$ it can be proved
in a similar way, using a block argument.

\subsection{Proof of Theorem \ref{thm:radius}}
\begin{proof}
As it is pointed out in \cite{BPSV}, it is very easy to see that 

\begin{equation} \label{thm1-lbd}
P(\hat{R}\geq n) \geq P_{cr}(0\leftrightarrow\partial B(n)),
\end{equation}
since
the whole $p_c$-open cluster of the origin is invaded before \emph{any} 
edge with $\tau$ value larger than $p_c$ is added to the invasion cluster.
To prove that the l.h.s. of \eqref{thm1-lbd} is smaller than some constant $c$ times the r.h.s.
is more involved. First note that it suffices to prove 
this for the case that $n$ is of the form $2^k$. Indeed, if it holds for those special cases then,
for any $2^{k-1}<n<2^k$ we have
\eqnsplst
{ P(\hat{R}\geq n)
  &\leq P(\hat{R}\geq 2^{k-1}) 
  \leq c P_{cr}\big[0 \cnctd \partial B(2^{k-1})\big] \\
  &\stackrel{\eqref{csakegyszerkell}}{\leq} \bar{c} P_{cr} \big[ 0 \cnctd \partial B(2^k) \big]
  \leq \bar{c} P_{cr} \big[0 \cnctd \partial B(n)\big]. }

First some additional notation and definitions. 
As in \cite{jarai} we define 
$\log^{(0)}k=k$ and $\log^{(j)}k=\log(\log^{(j-1)}k)$ for
all $j\geq 1$, as long as the right-hand side is well defined. For $k>10$
let
\begin{eqnarray}
\log^*k=\min\{j>0\ :\ \text{$\log^{(j)}k$ is well-defined and $\log^{(j)}k \leq 10$} \},
\end{eqnarray}
where the  choice of the constant 10 is quite arbitrary. Clearly,
$\log^{(j)}k>2$ for $j=0,1,\ldots,\log^{*}k$
and $k>10$. Further,
\begin{eqnarray}\label{pdef}
p_k(j):= \inf\Big\{p>p_c\ :\ L(p)\leq\frac{2^k}{C_2\log^{(j)}k}\Big\},
\end{eqnarray}
where the constant $C_2$ will be chosen later. 
It is easy to see that
$p_k(j)$ is well-defined for all sufficiently large $k$ (in fact, for all $k$ with 
$2^k > C_2 k$), and that the sequence $\{p_k(j)\}_{j=0}^{\log^*k}$ is decreasing in $j$.
The definition of $p_k(j)$ together with the right continuity of $L(p)$ 
and \eqref{smalljump} readily implies that
\begin{eqnarray}
\label{corineq}
C_2\log^{(j)}k\leq\frac{2^k}{L(p_k(j))}\leq DC_2\log^{(j)}k.
\end{eqnarray}
Now we decompose $P(\hat{R}\geq n)$ according to the value of $\hat{\tau}$ as follows, where we
note that since $\tau$ has a continuous distribution, $\hat{\tau}$
does not coincide with $p_k(j)$ for any $j=0, \ldots, \log^*k$, almost surely. 
\eqnspl{summa}
{ P(\hat{R} \geq n)
  &= P(\hat{R} \geq n, p_k(0)<\hat{\tau})
       + P(\hat{R}\geq n, \hat{\tau} < p_k(\log^*k)) \\
  &\quad+ \sum_{j=0}^{\log^*k-1} P(\hat{R}\geq n, p_k(j+1)<\hat{\tau}<p_k(j)). }
To bound the terms in \eqref{summa} we will use the following observations made in \cite{BPSV}.
Let $p$ be an arbitrary number between $p_c$ and 1. \\

\noindent
\emph{Observations}
\begin{itemize}
\item[(a)]  $\hat{\tau}<p$ if and only if the origin belongs to an infinite $p$--open cluster.
\item[(b)] If $\hat{\tau}>p$ and $\hat{R}\geq n$, then there is a $p$--closed circuit around $O$
in the dual lattice with diameter at least $n$.
\end{itemize}
The event in observation (b) will be denoted by $A_{n,p}$.
\begin{eqnarray*}
A_{n,p}:=\big\{ \text{$\exists$ $p$-closed circuit around $O$ in the dual with diameter at least $n$} \big\}.
\end{eqnarray*}
Starting with the first term of (\ref{summa}), Observation (b) gives 
\begin{eqnarray}\label{firsttermineq}
P(\hat{R}\geq n, p_k(0)<\hat{\tau})\leq P( A_{n,p_k(0)}).
\end{eqnarray}
It is well-known (see \cite{BPSV}
for more explanation and references) that there exist $C_3$ and $C_4$ such that
for all $p>p_c$,
\begin{eqnarray}\label{expdecay}
P(A_{n,p})\leq C_3\exp\Big\{-\frac{C_4n}{L(p)}\Big\}
\end{eqnarray}
Using the lower bound in (\ref{corineq}) and the definition of
$\log^{(0)}k$ we get that
\begin{eqnarray}\label{anpbound}
P(A_{n,p_k(0)})\leq
C_3\exp\Big\{-\frac{C_4n}{L(p_k(0))}\Big\}\stackrel{(\ref{corineq})}{\leq}C_3n^{-C_4 C_2}
\end{eqnarray}
As mentioned above, we have $P_{cr}(0\leftrightarrow \partial B(n))\geq C n^{-1/2}$. 
Hence, by taking $C_2 \geq 1/C_4$, we can ensure that 
\eqnst
{ P(A_{n,p_k(0)}) 
  \leq C_3 n^{-1} 
  \leq \tilde{C_3} P_{cr}(0 \shcnctd \partial B(n)). } 
{\bf Remark:} {\it For future purpose we will even take $C_2 \geq 2/C_4$}. 

\smallskip\noindent
For the second term of (\ref{summa}) we apply observation (a) to get
\begin{equation*}
P(\hat{R}\geq n, \hat{\tau}<p_k(\log^*k))\leq
P(\hat{\tau}<p_k(\log^*k)) \stackrel{\text{Obs. (a)}}{\leq} \theta(p_k(\log^*k)).
\end{equation*}
Furthermore, using (\ref{thetacorrineq}), \eqref{corineq}, the
definition of $p_k(\log^*k)$ and \eqref{csakegyszerkell},
we have
\eqnsplst
{ \theta(p_k(\log^*k)) 
  &\leq C_1 P_{cr}\big[0 \shcnctd \partial B(L(p_k(\log^*k))\big] \\
  &\leq C_1 P_{cr}\big[0 \shcnctd \partial B(\frac{2^k}{10 D C_2})\big] \\
  &\leq C_5 P_{cr}\big[0 \shcnctd \partial B(n)\big], }
for some constant $C_5$.

Now let us consider a typical term in the summation in
(\ref{summa}).
The two observations a few lines below \eqref{summa} (and the definition of $A_{n,p}$) give
\eqnspl{twoterms}
{ &P (\hat{R} \geq n, p_k(j+1)<\hat{\tau}<p_k(j)) \\ 
  &\qquad \leq  P (0 \pcnctd{p_k(j)} \infty,\, A_{n,p_k(j+1)} )\\
  &\qquad \leq \theta(p_k(j)) P(A_{n,p_k(j+1)}), }
where in the last inequality we use the Harris-FKG 
inequality \cite[Section 2.2]{grimmett}.
To bound the first factor in the right hand side of (\ref{twoterms}), note that
\eqnspl{sumterm1}
{ \theta(p_k(j))
  &\stackrel{(\ref{thetacorrineq})}{\leq} C_1P_{cr}(0\leftrightarrow L(p_k(j)))\\
  &= C_1P_{cr}(0\leftrightarrow \partial B(2^k))
      \frac{P_{cr}(0\leftrightarrow L(p_k(j)))}{P_{cr}(0\leftrightarrow \partial B(2^k))}\\
  &\stackrel{(\ref{csakegyszerkell})}{\leq} \frac{C_1}{D_1}
      P_{cr}(0\leftrightarrow \partial B(2^k)) \bigg(\frac{2^k}{L(p_k(j))}\bigg)^{1/2}\\
  &\stackrel{(\ref{corineq})}{\leq} \frac{C_1}{D_1}
      P_{cr}(0\leftrightarrow \partial B(2^k))(DC_2\log^{(j)}k)^{1/2}. }
The second factor in the right hand side of (\ref{twoterms}) can be bounded using (\ref{expdecay}),
\eqref{corineq}, \eqref{pdef} and the choice of $C_2$:
\begin{eqnarray}\label{sumterm2}
P(A_{n,p_k(j+1)})\leq
C_3\exp\Big\{-\frac{C_4n}{L(p_k(j+1))}\Big\}\leq
C_3(\log^{(j)}k)^{-1},
\end{eqnarray}
Combining (\ref{sumterm1}) and (\ref{sumterm2}) gives
\begin{eqnarray}
\theta(p_k(j))P(A_{n,p_k(j+1)})\leq
C_8(\log^{(j)}k)^{-1/2}P_{cr}\big[0\leftrightarrow\partial
B(n)\big].
\end{eqnarray}
To conclude the proof it suffices to show that
\begin{eqnarray}\label{supsum}
\sup_{k>10}\sum_{j=0}^{\log^*k-1}(\log^{(j)}k)^{-1/2} < \infty.
\end{eqnarray}
Recall from the definitions that $\log^{(j)}k>2$. Applying this to the case $j=\log^*k$ shows that the
last term in the sum in (\ref{supsum}) is at most $(e^2)^{-1/2}$. Similarly, the penultimate term is
at most $(\exp(e^2))^{-1/2}$, etc. This leads to the finite upper bound
$C_9 := \frac{1}{\sqrt{e^2}}+\frac{1}{\sqrt{e^{e^2}}}+\ldots$ for the l.h.s. of \eqref{supsum}. \\
Putting everything together we get
\begin{eqnarray*}
P(\hat{R}\geq n)\leq
\big(\tilde{C}_3+C_5+C_8C_9\big)P_{cr}\big[0\leftrightarrow\partial
B(n)\big]. 
\end{eqnarray*}\end{proof}

\subsection{Proof of Theorem \ref{thm:volume}}

For short, we use the following notation:
\eqnsplst
{ \pi(n) 
  &= P_{cr}(0 \shcnctd \partial B(n)); \\
   \pi(n,p) 
   &= P_p ( 0 \shcnctd \partial B(n)). }
Recall that $s(n) = n^2 \pi(n)$.

The difficult part of Theorem \ref{thm:volume} is the third inequality. We need the following
key ingredient. 

\begin{lemma}
\label{lem:largedev}
There exist constants $C_{10}$ and $C_{11}$, such that 
\eqnsplst
{ &P_p \left( 0 \shcnctd \infty,\, | \cC(0) \cap B(2^k) | > s(n) \right) \\
  &\qquad \le \theta(p)\, 2 C_{10}\, 
     \exp \left\{ - (2 C_{11})^{-1} \frac{s(n)}{2^{2k} \pi(2^k, p)} \right\}, \quad
     p > p_c, 2^k \le n. }
\end{lemma}

\begin{proof}
The proof is based on the following moment estimate:
\eqn{e:moments}
{ E_p \Big( |\cC(0) \cap B(2^k)|^t  \Big| 0 \cnctd \infty \Big)
   \le C_{10}\, t!\, \left[ C_{11}\, 2^{2k}\, \pi(2^k, p) \right]^t, \quad t \ge 1. }
Very similar estimates were proved in \cite[Theorem (8)]{kesten86} and in \cite{nguyen87}.
To adapt their proofs in order to obtain \eqref{e:moments}, one
merely needs that the inequality $\sum_{m = 0}^n \pi(n,p) \le C n \pi(n,p)$ holds
for all $p \ge p_c$ (with some constant $C$ independent of $p$).
From \eqref{e:moments}, we readily get
\eqnst
{ E_p \Bigg( \exp \left\{ \lambda \frac{ | \cC(0) \cap B(2^k) | }{ 2^{2k}\, \pi(2^k, p) } \right\} \,\Bigg|\,
           0 \cnctd \infty \Bigg)
   \le C_{10} \frac{1}{1 - \lambda C_{11}},\quad 0 < \lambda < C_{11}^{-1}.  } 
Taking $\lambda = (2 C_{11})^{-1}$ we easily obtain the estimate of the lemma.   
\end{proof}

\begin{proof}[Proof of Theorem \ref{thm:volume}]
The first inequality follows from \cite[Remark (9)]{kesten86}.
The second inequality follows immediately from the fact that the 
$p_c$-open cluster containing the origin is a subset of $\hat{V}$.

The third inequality will be proved by a decomposition, somewhat similar to 
the one in Theorem \ref{thm:radius}, but now two-fold: this time we will also decompose
according to the value of $\hat{R}$. As in the proof of 
Theorem \ref{thm:radius}, without loss of generality we may assume 
that $n$ is of the form $2^N$. 

Let
\eqnst
{ E_{n,k}
   = \{ 2^{k-1} < \hR \le 2^k,\, |\hV| > s(n) \}. }
Note that $s(n) \ge C_{12} n^{3/2}$, and 
$|B(2^k)| \le C_{13} 2^{2k}$. Letting  
\eqnst
{ k_0 
  := \max \{k : C_{13} 2^{2k} \le C_{12} n^{3/2} \}, }
for $k < k_0$, $\hR \le 2^k$ implies $|\hV| \le C_{13} 2^{2k} \le s(n)$, and hence 
$E_{n,k} = \emptyset$. Therefore, we can write
\eqnspl{e:decomp1}
{ P( |\hat{V}| > s(n) ) 
  &\le P( \hat{R} > n )
     + \sum_{k = k_0}^{N}
     P( E_{n,k} ). }
The first term on the right hand side is at most $C_{14} \pi(n)$, by
Theorem \ref{thm:radius}. Consider now a general term of the sum. 
We decompose this according to the value of
$\hat\tau$ as follows:
\eqnspl{e:decomp2}
{ &P( 2^{k-1} < \hat{R} \le 2^k,\, |\hat{V}| > s(n) ) \\
  &\qquad = P( E_{n,k},\, \hat\tau > p_k(0) ) 
    + \sum_{j = 0}^{\log^* k}
    P( E_{n,k},\, 
    p_k(j+1) < \hat\tau < p_k(j) ), }
where we let $p_k(\log^* k + 1) = p_c$.

We first look at the event in the first term on the right hand side. This event implies
the occurrence of $A_{ 2^{k-1}, p_k(0)}$. Hence, by virtue of  \eqref{anpbound},
its probability is at most $C_{15} (2^{2k})^{-\tilde{C_4} C_2}$. 
By the choice of $C_2$, we have
$\tilde{C_4} C_2 \ge 1$. Hence the sum over $k_0 \le k \le N$ is bounded by 
$C_{16} (2^{2 k_0})^{-1}$. By the definition of $k_0$, this is $o(\pi(n))$.

Consider now the event in the general term on the right hand side of 
\eqref{e:decomp2}. This event implies the following two events:
\begin{itemize}
\item[(i)] $A_{2^{k-1}, p_k(j+1)}$;
\item[(ii)] $\{ 0 \stackrel{p_k(j)} \longleftrightarrow \infty,\, 
  |\cC(0; p_k(j)) \cap B(2^k)| > s(n) \}$;
\end{itemize}
where $\cC(0 ; p)$ denotes the $p$-open cluster of $0$.
Since (i) is a decreasing  and (ii) an increasing event,
the Harris-FKG inequality yields 
that the general term in \eqref{e:decomp2} is at most the product of the probabilities of
event (i) and event (ii). \\
As to event (i), the same arguments that led to \eqref{sumterm2} (and noting the Remark a few lines
below \eqref{anpbound}) show that for $j < \log^* k$ this has probability less than or equal to 

\begin{equation}
\label{boundi}
C_3 (\log^{(j)} k)^{-1}
\end{equation}
It is easy to see that, after increasing the value of $C_3$ if necessary, this bound even holds for
$j = \log^* k$.

As to event (ii), by Lemma \ref{lem:largedev} this has probability at most
\begin{equation}
\label{boundiia}
\theta(p_k(j)) (2 C_{10}) \exp\left\{ - (2 C_{11})^{-1} \frac{ s(n) }{ 2^{2k}\, \pi(2^k, p_k(j)) } \right\}.
\end{equation}
Applying the first inequality in \eqref{thetacorrineq} to the probability in the exponent in
\eqref{boundiia}, and then applying \eqref{sumterm1} twice, shows that \eqref{boundiia} is at most a 
constant times
\begin{equation}
\label{boundiib}
\pi(2^k) (\log^{(j)} k)^{1/2} \exp \left\{ - C_{18} \frac{ 2^{2N} \pi(n) }{ 2^{2k} \pi(2^k) }
 (\log^{(j)} k)^{-1/2}  \right\}.
\end{equation}
Combining this with \eqref{boundi} gives that the general term in \eqref{e:decomp2} is at most
\begin{equation}
\label{boundcomb}
C_{17} \pi(n) (\log^{(j)} k)^{-1/2} \frac{\pi(2^k)}{\pi(n)}
     \exp \left\{ - C_{18} \frac{ 2^{2N} \pi(n) }{ 2^{2k} \pi(2^k) }  (\log^{(j)} k)^{-1/2}  \right\}.
\end{equation}
Due to \eqref{csakegyszerkell}, this as at most
\eqn{e:tosum}
{ C_{19} \pi(n) (\log^{(j)} k)^{-1/2} 2^{(N-k)/2}
      \exp \left\{ - C_{20} 2^{(N - k)(3/2)} (\log^{(j)} k)^{-1/2}  \right\}. }

We split the sums over $j$ and $k$ into two parts:
\begin{itemize}
\item[(1)] $2^{(N - k)} \le (\log^{(j)} k)^{1/2}$;
\item[(2)] $2^{(N - k)} > (\log^{(j)} k)^{1/2}$.
\end{itemize}
In case (1), we bound the exponential
in \eqref{e:tosum} by $1$, and we have
\eqnst
 { (\log^{(j)} k)^{-1/2} 2^{(N - k)/2} \le (\log^{(j)} k)^{-1/4} \le C_{21} (\log^{(j)} N)^{-1/4}.}
The number of possible values of $k$ is at most
\eqnst
{ (2 \log 2)^{-1} \log^{(j+1)} k \le C_{22} (\log^{(j)} N)^{1/8}. } 
Hence the contribution of this case is bounded by
\eqnst
{ \sum_{j = 0}^{\log^* N} ( \log^{(j)} N )^{-1/8} 
  \le C_{23}. }
In case (2), we bound the exponential by $\exp \{ - C_{20} 2^{(N - k)/2} \}$, and we have
$(\log^{(j)} k)^{-1/2} 2^{(N - k)/2} \le 2^{(N - k)/2}$. The sum over $k$ can be bounded 
as follows:
\eqnst
{ \sum_{k : N - k \ge c \log^{(j+1)} N} 2^{(N - k)/2} \exp \{ - C_{20} 2^{(N - k)/2} \}
  \le C_{24} \exp \{ - C_{25} (\log^{(j)} N)^{c_1} \}, }
for some $c_1 > 0$. The sum of the right hand side over $j$ is again bounded.
This proves the theorem.
\end{proof}


\begin{thebibliography}{00}

\bibitem{AGdHS06}Angel, O, Goodman, J, den Hollander, F. and Slade, G.:
Invasion percolation on regular trees, preprint, arXiv:math/0608132v1, (2006).

\bibitem{BK}van den Berg, J., Kesten, H.: Inequalities with application 
to percolation and reliability.
J. Appl. Prob. \textbf{22}, 556--569 (1985)

\bibitem{BPSV} van den Berg, J., Peres, Y., Sidoravicius, V. and Vares, M.E.:
Random spatial growth with paralyzing obstacles, preprint, arXiv:0706.0219 (2007).

\bibitem{newman} Chayes, J.T., Chayes, L. and Newman, C.M.: Stochastic
geometry of invasion percolation, Commun. Math. Phys. \textbf{101},
383--407 (1985)

\bibitem{newman2} Chayes, J.T., Chayes, L. and Newman, C.M.:
Bernoulli percolation above threshold: An invasion percolation
analysis. Ann. Probab. \textbf{15}, 1272--1287 (1987)

\bibitem{grimmett} Grimmett, G.R.: \textit{Percolation},
2nd edition. Springer--Verlag (1999)

\bibitem{jarai} J\'{a}rai, A.A.: Invasion percolation and the
incipient infinite cluster in 2D. Commun. Math. Phys. \textbf{236},
311--314 (2003)

\bibitem{kesten86} Kesten, H.: The incipient infinite cluster in 
two-dimensional percolation.  
Probab. Theory Related Fields  {\bf 73}, 369--394 (1986)

\bibitem{kesten} Kesten, H.: Scaling relations for 2D percolation.
Commun. Math. Phys. \textbf{109}, 109--156 (1987)

\bibitem{LPS06} Lyons, R., Peres, Y.~and Schramm, O.: Minimal spanning 
forests. Ann. Probab. \textbf{34}, 1665--1692 (2006)

\bibitem{nguyen87} Nguyen, B.G.:
Typical cluster size for two-dimensional percolation processes.
J. Stat. Phys. \textbf{50}, 715--726 (1987)

\bibitem{LSW02} Lawler, G.F., Schramm, O. and Werner, W.:
One-arm exponent for critical 2D percolation.
Electron. J. Probab. \textbf{2}, 13 pp., electronic (2002)

\bibitem{St-New} D.L. Stein and C.M. Newman. Broken ergodicity and the geometry of rugged landscapes,
{\it Phys. Rev. E} {\bf 51}, 5228--5238 (1995)

\bibitem{wilkinson} Wilkinson, D. and Willemsen, J.F.: Invasion
percolation: A new form of percolation theory. J. Phys. A.
\textbf{16}, 3365--3376 (1983)

\end{thebibliography}
\end{document}